\documentclass{ert-l}
\usepackage{amsmath, amssymb, amscd}

\newtheorem{thm}{Theorem}

\newtheorem{rem}{Remark}

\newtheorem{df}{Definition}
\newtheorem{ex}[df]{Example}

\renewcommand{\phi}{\varphi}
\renewcommand{\epsilon}{\varepsilon}

\newcommand{\cal}{\mathcal}
\newcommand{\BB}{\mathbb}
\newcommand{\g}{\mathfrak}
\newcommand{\separate}{\vskip5pt}

\newcommand{\supp}{\operatorname{supp}}

\newcommand{\re}{\operatorname{Re}}

\newcommand{\End}{\operatorname{End}}
\newcommand{\rhom}{\operatorname{RHom}}
\newcommand{\modg}{\operatorname{Mod}_{quasi-G_{\BB C}-eq}}
\newcommand{\modge}{\operatorname{Mod}_{G_{\BB C}-eq}}
\newcommand{\thom}{\operatorname{Thom}}
\newcommand{\parag}{${\mathcal x}$}

\begin{document}

\title[Riemann-Roch-Hirzebruch integral formula]
{Riemann-Roch-Hirzebruch integral formula for characters of
reductive Lie groups}

\author{Matvei Libine}

\address{Department of Mathematics and Statistics, University of Massachusetts,
Lederle Graduate Research Tower, 710 North Pleasant Street, Amherst, MA 01003}

\email{matvei@math.umass.edu}



\begin{abstract}
Let $G_{\BB R}$ be a real reductive Lie group acting on a manifold $M$.
M.~Kashiwara and W.~Schmid in \cite{KaSchm} constructed representations
of $G_{\BB R}$ using sheaves and quasi-$G_{\BB R}$-equivariant
${\cal D}$-modules on $M$.
In this article we prove an integral  character formula for these
representations (Theorem \ref{main}).
Our main tools will be the integral localization formula recently proved
in \cite{L3} and the integral character formula proved by W.~Schmid and
K.~Vilonen in \cite{SchV2} (originally established by W.~Rossmann in \cite{Ro})
in the important special case when the manifold $M$ is the flag variety of
$\BB C \otimes_{\BB R} \g g_{\BB R}$ -- the complexified Lie algebra of
$G_{\BB R}$.
In the special case when $G_{\BB R}$ is commutative and the
${\cal D}$-module is the sheaf of sections of a
$G_{\BB R}$-equivariant line bundle over $M$ this integral character
formula will reduce to the classical Riemann-Roch-Hirzebruch formula.
As an illustration we give a concrete example on the enhanced flag variety.
\end{abstract}

\maketitle

\tableofcontents

\begin{section}
{Introduction}
\end{section}

We consider pairs of Lie groups: a real linear reductive group $G_{\mathbb R}$
sitting inside its complexification $G_{\mathbb C}$.
For example:
$$
\begin{matrix}
GL(n, \mathbb R) & \subset & GL(n, \mathbb C) \\
GL^+(n, \mathbb R) & \subset & GL(n, \mathbb C) \\
U(n) & \subset & GL(n, \mathbb C)
\end{matrix}
\qquad
\begin{matrix}
SL(n, \mathbb R) & \subset & SL(n, \mathbb C) \\
SO(n, \mathbb C) & \subset & SL(n, \mathbb C)  \\
SU(n) & \subset & SL(n, \mathbb C)
\end{matrix}
\qquad
Sp(n, \mathbb R) \subset Sp(n, \mathbb C).
$$
Suppose that $G_{\mathbb C}$ acts algebraically on a smooth complex
projective variety $M$.
Fix a $G_{\mathbb C}$-invariant open algebraic subset $U \subset M$,
and take a $G_{\mathbb C}$-equivariant algebraic line bundle
$({\bf E}, \nabla_{\bf E})$ over $U$ with a $G_{\mathbb C}$-invariant
algebraic flat connection $\nabla_{\bf E}$.
Let $S_{\mathbb R} \subset M$ be an open $G_{\mathbb R}$-invariant subset
(which may or may not be $G_{\mathbb C}$-invariant) and
consider the cohomology spaces
\begin{equation}  \label{hspaces}
H^p(S_{\mathbb R}, {\cal O}({\bf E})), \qquad p \in \mathbb Z.
\end{equation}

The classical Riemann-Roch-Hirzebruch formula computes the index of
${\bf E}$, i.e. the alternating sum
$$
\sum_p (-1)^p \dim H^p(S_{\mathbb R}, {\cal O}({\bf E}))
$$
with $S_{\mathbb R} = U = M$.
For general $S_{\mathbb R}$ and $U$, however, these dimensions can be infinite.
To work around this problem we regard the vector spaces (\ref{hspaces})
as representations of $G_{\mathbb R}$, and, as a substitute for the index,
we ask for the character of the virtual representation
$$
\sum_p (-1)^p H^p(S_{\mathbb R}, {\cal O}({\bf E})).
$$
(Recall that for finite-dimensional representations the value of the character
at the identity element $e \in G$ equals the dimension of the representation.)

In this article we establish an integral formula for characters of
virtual representations of this kind in a more general setting
of ${\cal D}$-modules.
Existence of such integral character formula was conjectured by W.~Schmid.
My work toward this formula has led to the integral localization formula
proved in \cite{L3}.

We use the following convention: whenever $A$ is a subset of $B$,
we denote the inclusion map $A \hookrightarrow B$ by $j_{A \hookrightarrow B}$.

\separate

\begin{section}
{Preliminary results}
\end{section}

Let $G_{\BB C}$ be a connected complex algebraic reductive group
which is defined over $\BB R$.
We will be primarily interested in representations of a subgroup
$G_{\BB R}$ of $G_{\BB C}$ lying between the group of real points
$G_{\BB C}(\BB R)$ and the identity component $G_{\BB C}(\BB R)^0$.
We regard $G_{\BB R}$ as a real reductive Lie group.
Let $\g g_{\BB C}$ and $\g g_{\BB R}$ be the respective Lie algebras
of $G_{\BB C}$ and $G_{\BB R}$.
Let ${\cal B}$ denote the flag variety of $\g g_{\mathbb C}$.
It is a smooth complex projective variety consisting of all Borel
subalgebras of $\g g_{\BB C}$. The group $G_{\mathbb C}$ acts on ${\cal B}$
transitively.

If the group $G_{\mathbb R}$ is compact, then all irreducible representations
of $G_{\mathbb R}$ can be enumerated by their highest weights $\lambda$
lying in the weight lattice $\Lambda \subset i \g g_{\mathbb R}^*$.
The Borel-Weil-Bott theorem can be regarded as an explicit
construction of a holomorphic $G_{\mathbb R}$-equivariant line bundle
${\cal L}_{\lambda} \to {\cal B}$ such that the resulting
representation of $G_{\mathbb R}$ in the cohomology groups is:
\begin{eqnarray*}
&& H^p({\cal B}, {\cal O}({\cal L}_{\lambda}))=0 \text{\quad if $p \ne 0$}, \\
&& H^0({\cal B}, {\cal O}({\cal L}_{\lambda})) \simeq \pi_{\lambda},
\end{eqnarray*}
where ${\cal O}({\cal L}_{\lambda})$ is the sheaf of sections of
${\cal L}_{\lambda}$ and
$\pi_{\lambda}$ denotes the irreducible representation of
$G_{\mathbb R}$ of highest weight $\lambda$.
Then N.~Berline and M.~Vergne \cite{BV}, \cite{BGV} observed that the
character of $\pi_{\lambda}$, as a function on $\g g_{\mathbb R}$,
can be expressed as an integral over ${\cal B}$ of a certain naturally
defined equivariantly closed form.
They proved it by applying their famous integral localization
formula and matching contributions from fixed points with terms
of the Weyl character formula.
(This is a restatement of Kirillov's character formula.)

\separate

M.~Kashiwara and W.~Schmid \cite{KaSchm}
generalize the Borel-Weil-Bott construction.
Instead of line bundles over the flag variety ${\cal B}$ they consider
$G_{\mathbb R}$-equivariant sheaves\footnote{Strictly speaking,
${\cal F}$ is not a sheaf on ${\cal B}$
but rather an element of the ``$G_{\mathbb R}$ equivariant derived category
on ${\cal B}$ with twist $(-\lambda-\rho)$'' denoted by
$\operatorname{D}_{G_{\BB R}} ({\cal B})_{-\lambda}$
(see \cite{SchV2} and Remark \ref{twistings}).}
${\cal F}$ and, for each integer
$p \in \BB Z$, they define representations of $G_{\mathbb R}$ in
$\operatorname{Ext}^p({\cal F},{\cal O}_{{\cal B}^{an}})$, where
${\cal O}_{{\cal B}^{an}}$ denotes the sheaf of analytic functions on
${\cal B}$.
In other words, M.~Kashiwara and W.~Schmid prove that these vector spaces
possess a natural Fr\'echet topology and the resulting representations of
$G_{\mathbb R}$ are admissible of finite length.
Hence these representations have characters in the sense of Harish-Chandra.
Let $\tilde \theta$ be the character of the virtual representation of
$G_{\mathbb R}$
$$
\sum_p (-1)^p \operatorname{Ext}^p
\bigl( \BB D_{\cal B} {\cal F}, {\cal O}_{{\cal B}^{an}}(\lambda) \bigr),
$$
where $\BB D_{\cal B} {\cal F}$ denotes the Verdier dual of ${\cal F}$,
${\cal O}_{{\cal B}^{an}}$ 
and $\lambda$ is some twisting parameter lying in $\g h_{\BB C}^*$ --
the dual space of the universal Cartan algebra of $\g g_{\BB C}$.
Here the character $\tilde \theta$ is a distribution on $\g g_{\mathbb R}$.
Then W.~Schmid and K.~Vilonen \cite{SchV2} prove
two character formulas for $\tilde \theta$.

The integral character formula expresses $\tilde \theta$ as an integral of a
certain differential form (independent of ${\cal F}$) over the characteristic
cycle $Ch({\cal F})$ of ${\cal F}$.
Characteristic cycles were introduced by M.~Kashiwara and their
definition can be found in \cite{KaScha}. On the other hand,
W.~Schmid and K.~Vilonen give a geometric way to understand
characteristic cycles in \cite{SchV1}.
A comprehensive treatment of characteristic cycles can be found in \cite{Schu}.
The cycle $Ch({\cal F})$ is a conic Lagrangian Borel-Moore homology cycle
lying inside the cotangent space $T^*{\cal B}$.
If the sheaf ${\cal F}$ happens to be perverse,
the characteristic cycle of ${\cal F}$ equals the characteristic cycle of the
holonomic ${\cal D}$-module corresponding to ${\cal F}$ via the
Riemann-Hilbert correspondence.
Originally existence of such
formula over an unspecified cycle was established by W.~Rossmann in \cite{Ro}.
In the special case when the group $G_{\mathbb R}$ is compact, the integral
character formula reduces to Kirillov's character formula.

On the other hand, Harish-Chandra showed that the distribution $\tilde\theta$
is given by integration against a certain function $F_{\tilde\theta}$
on $\g g_{\mathbb R}$.
This function $F_{\tilde\theta}$ is an $Ad(G_{\mathbb R})$-invariant locally
$L^1$-function, its restriction to the set of regular semisimple
elements $\g g_{\mathbb R}^{rs}$ of $\g g_{\mathbb R}$ can be represented
by an analytic function.
According to Harish-Chandra, $F_{\tilde\theta}$ can be expressed as follows.
For a regular semisimple element $X \in \g g_{\mathbb R}^{rs}$ we denote by
$\g t_{\mathbb C}(X) \subset \g g_{\mathbb C}$ the unique (complex) Cartan
algebra containing $X$ and by $\Psi(X) \subset \g t_{\mathbb C}(X)^*$ we
denote the root system of $\g g_{\mathbb C}$ with respect to
$\g t_{\mathbb C}(X)$.
Since $X$ is regular semisimple, the number of Borel subalgebras
$\{\g b_1,\dots,\g b_{|W|}\}$ containing $X$ is exactly the order of the
Weyl group $W$.
Each Borel subalgebra $\g b_k \in {\cal B}$ containing $X$
(and hence containing $\g t_{\mathbb C}(X)$) determines a positive root system 
$\Psi_{\g b_k}^+(X) \subset \Psi(X)$ consisting of all those roots whose root
spaces are {\em not} contained in $\g b_k$;
that is $\g b_k$ contains all negative root spaces relative to
$\Psi_{\g b_k}^+(X)$. Then
\begin{equation}  \label{fpcharformula}
F_{\tilde\theta} (X) =
\sum_{k=1}^{|W|}
m_{\g b_k}(X) \cdot \frac {e^{\langle X,\lambda_{\g b_k} \rangle}}
{\prod_{\alpha \in \Psi_{\g b_k}^+(X)} \alpha(X)},
\end{equation}
where $m_{\g b_k}(X)$'s are some integer multiplicities.
W.~Schmid and K.~Vilonen give a formula for $m_{\g b_k}(X)$'s in terms of
local cohomology of ${\cal F}$ (formula (5.25b) in \cite{SchV2}).
Pick another positive root system $\Psi^{\le}(X) \subset \Psi(X)$ such that
$$
\re \alpha (X) \le 0 \qquad \text{for all $\alpha \in \Psi^{\le}(X)$.}
$$
(If there are roots in $\Psi$ which take purely imaginary values on $X$,
there will be different possible choices for $\Psi^{\le}(X)$.)
Let $B(X) \subset G_{\BB C}$ be the Borel subgroup whose Lie algebra
contains $\g t_{\BB C}(X)$ and the positive root spaces corresponding
to $\Psi^{\le}(X)$.
The action of $B(X)$ on ${\cal B}$ has precisely $|W|$ orbits
$O_1, \dots, O_{|W|}$. Each of these orbits $O_k$ is a locally closed
subset of ${\cal B}$ and contains exactly one of the Borel subalgebras
$\{\g b_1,\dots,\g b_{|W|}\}$ containing $X$.
We order the Borel subalgebras containing $X$ and the orbits of $B(X)$ so that
$\g b_k$ is contained in $O_k$, $k=1,\dots,|W|$.
Then
\begin{equation}  \label{m}
m_{\g b_k}(X) = \chi
\bigl( {\cal H}^*_{O_k} (\mathbb D_{\cal B} {\cal F})_{\g b_k} \bigr)
= \chi \bigl( (j_{\{\g b_k\} \hookrightarrow O_k})^* \circ
(j_{O_k \hookrightarrow {\cal B}})^! \mathbb D_{\cal B} {\cal F}) \bigr)
\end{equation}
(recall that $\BB D_{\cal B} {\cal F}$ denotes the Verdier dual of ${\cal F}$).
These multiplicities are exactly the local contributions of points
$\g b_k \in {\cal B}$ to the Lefschetz fixed point formula, as generalized
to sheaf cohomology by M.~Goresky and R.~MacPherson \cite{GM}.
W.~Schmid and K.~Vilonen call (\ref{fpcharformula}) combined with (\ref{m})
the fixed point character formula because the set of all Borel subalgebras
$\{\g b_1,\dots,\g b_{|W|}\}$ containing $X$ can be expressed as the set of
zeroes of the vector field generated by the infinitesimal action of
$X$ on ${\cal B}$.
In the special case when the group $G_{\mathbb R}$ is compact,
and  $\lambda+\rho$ is an integral weight, all the multiplicities
$m_{\g b_k}(X)$ are equal to each other, say,
$$
m_{\g b_1}(X) = \dots = m_{\g b_{|W|}}(X) = \kappa \in \mathbb Z,
$$
and (\ref{fpcharformula}) reduces to $\kappa$ times the Weyl character formula.
The fixed point formula was conjectured by M.~Kashiwara
\cite{Ka}, and its proof uses the above-mentioned generalization of the
Lefschetz fixed point formula to sheaf cohomology \cite{GM}.

There is a striking relationship between these two character formulas.
In the compact group case N.~Berline and M.~Vergne \cite{BV}, \cite{BGV}
gave a simple proof of Kirillov's character formula using their
integral localization formula; they matched contributions from zeroes of
vector fields with terms of the Weyl character formula.
However, in the non-compact group case their argument breaks down because
their localization formula works for compact groups only.
Originally W.~Schmid and K.~Vilonen \cite{SchV2} proved these
character formulas independently of each other using representation
theory methods. Thus, besides an important representation-theoretical result,
they formally established existence of an integral localization formula
in a very special case.
In the announcement \cite{Sch} W.~Schmid posed a question:
``Can this equivalence of character formulas be seen directly without
a detour to representation theory, just as in the compact case.''
In \cite{L1}, \cite{L2} I provide such a geometric link, then in \cite{L3}
I establish a general integral localization formula for non-compact group
actions.
In turn, this article uses the new localization formula for
non-compact group actions to give a generalization of the integral character
formula to representations associated to sheaves on manifolds other than the
flag manifold of $\g g_{\mathbb C}$.

\separate

We describe the integral character formula for $\tilde \theta$ in more detail;
its ingredients will be used in our character formula (\ref{mainintegral}).
The character $\tilde\theta$ is a distribution on
$\Omega_c^{top}(\g g_{\BB R})$ -- the space of smooth compactly
supported differential forms on $\g g_{\BB R}$ of top degree.
If $\phi \in \Omega_c^{top}(\g g_{\BB R})$, then we define its Fourier
transform $\hat \phi \in {\cal C}^{\infty}(\g g_{\BB C}^*)$
as in \cite{L1}, \cite{L2}, \cite{L3}, \cite{Ro} and \cite{SchV2}:
\begin{equation}  \label{FT}
\hat \phi (\xi) = \int_{\g g_{\BB R}}
e^{\langle X, \xi \rangle} \phi(X),
\qquad X \in \g g_{\BB R}, \: \xi \in \g g_{\BB C}^*,
\end{equation}
without the customary factor of $i = \sqrt{-1}$ in the exponent.

For a smooth complex algebraic manifold $M$ we denote by $\sigma_M$ the
canonical complex symplectic form on the holomorphic cotangent space $T^*M$.
In general, when the group $G_{\mathbb C}$ acts on a complex manifold $M$,
for $X \in \g g_{\BB C}$, we denote by $X_M$ the vector field on
$M$ given by (notice the minus sign)
\begin{equation}  \label{vfield}
(X_M \cdot \phi) (x) =
\frac d{d\epsilon} \phi \bigl( \exp (-\epsilon X)x \bigr) \Bigr|_{\epsilon=0},
\qquad \phi \in {\cal C}^{\infty}(M).
\end{equation}
Then the {\em moment map} on the holomorphic cotangent space
$\mu_M : T^*M \to \g g_{\BB C}^*$ is defined by
\begin{equation}  \label{mu}
\mu_M (\zeta): X \mapsto - \langle \zeta, X_M \rangle,
\qquad X \in \g g_{\BB C},\: \zeta \in T^*M.
\end{equation}
When the manifold $M$ is the flag variety ${\cal B}$ we get the complex
algebraic symplectic form $\sigma_{\cal B}$ on $T^*{\cal B}$, the moment map
$\mu_{\cal B}: T^*{\cal B} \to \g g_{\BB C}^*$ and the vector field
$X_{\cal B}$ for each element $X \in \g g_{\mathbb C}$.

We fix a compact real form (i.e. a maximal compact subgroup)
$U_{\BB R} \subset G_{\BB C}$ with Lie algebra $\g u_{\mathbb R}$.
We will use Rossmann's map
${\cal B} \ni x \mapsto \lambda_x \in \g g_{\BB C}^*$
defined in \cite{Ro} and Section 8 of \cite{SchV1}
(here we use the notations of \cite{SchV1} and \cite{SchV2}).
Recall that the twisting parameter $\lambda$ is an element of the dual of
the universal Cartan algebra $\g h_{\mathbb C}$ of $\g g_{\mathbb C}$; 
$\g h_{\mathbb C}$ is not a subalgebra of $\g g_{\mathbb C}$, but is
canonically isomorphic to any Cartan subalgebra
$\g t_{\mathbb C} \subset \g g_{\mathbb C}$ equipped with a specified choice
of positive root system $\Phi^+$.
For $x \in {\cal B}$, we denote by $\g b_x \subset \g g_{\mathbb C}$ the Borel
subalgebra corresponding to $x$. Then $\g h_{\mathbb C}$ is canonically
isomorphic to the quotient $\g b_x / [\g b_x, \g b_x]$, so that $\g b_x$
contains all the negative root spaces.
Thus we have an exact sequence of vector spaces
\begin{equation}  \label{split}
0 \to [\g b_x, \g b_x] \to \g b_x \to \g h_{\mathbb C} \to 0.
\end{equation}
In general, this sequence does not have a canonical splitting.
But once a choice of a compact real form $U_{\mathbb R} \subset G_{\BB C}$
is made, we can split this sequence as follows.
Let $T_{\mathbb R}(x)$ be the stabilizer of $x$ in $U_{\mathbb R}$, it is
a maximal torus in $U_{\mathbb R}$, and set
$$
\g t_{\mathbb R}(x) = Lie (T_{\mathbb R}(x)) \subset \g u_{\mathbb R} \subset
\g g_{\mathbb C},
\qquad
\g t_{\mathbb C}(x) = \g t_{\mathbb R}(x) \otimes_{\mathbb R} \mathbb C
\subset \g g_{\mathbb C}.
$$
Then $\g t_{\mathbb C}(x)$ is a Cartan subalgebra of $\g g_{\mathbb C}$
which lies in $\g b_x$ and comes equipped with a system of positive roots
$\Psi^+_{\g b_x}$ so that $\g b_x$ contains all the negative root spaces.
Hence we get a composition of maps
$$
\g h_{\mathbb C} \tilde \rightarrow \g t_{\mathbb C}(x) \hookrightarrow 
\g b_x
$$
which splits (\ref{split}):
$$
\g b_x \simeq \g t_{\mathbb C}(x) \oplus [\g b_x, \g b_x]
\simeq \g h_{\mathbb C} \oplus [\g b_x, \g b_x].
$$
We also get a splitting of $\g g_{\mathbb C}$:
$$
\g g_{\mathbb C} \simeq \g t_{\mathbb C}(x) \oplus
\Bigl(\bigoplus \text{root spaces of $\g t_{\mathbb C}(x)$} \Bigr)
\simeq \g h_{\mathbb C} \oplus
\Bigl(\bigoplus \text{root spaces of $\g t_{\mathbb C}(x)$} \Bigr).
$$
Taking duals we get a splitting of $\g g_{\mathbb C}^*$:
$$
\g g_{\mathbb C}^* \simeq \g t_{\mathbb C}(x)^* \oplus
\Bigl(\bigoplus \text{root spaces of $\g t_{\mathbb C}(x)$} \Bigr)^*
\simeq \g h_{\mathbb C}^* \oplus
\Bigl(\bigoplus \text{root spaces of $\g t_{\mathbb C}(x)$} \Bigr)^*.
$$
Hence $\lambda \in \g h_{\mathbb C}^*$ gets identified via this splitting
with an element $\lambda_x \in \g g_{\mathbb C}^*$.
The map $\lambda_x: {\cal B} \to \g g_{\BB C}^*$ is smooth, real algebraic
and $U_{\mathbb R}$-equivariant, but in general not $G_{\mathbb C}$- or
$G_{\mathbb R}$-equivariant, nor complex algebraic.
When $\lambda$ is regular, the {\em twisted moment map}
$$
\mu_{\lambda} =_{def} \mu_{\cal B} + \lambda_x :
T^*{\cal B} \to \g g_{\mathbb C}^*
$$
is a real algebraic $U_{\mathbb R}$-equivariant diffeomorphism of $T^*{\cal B}$
onto the coadjoint orbit of $G_{\mathbb C}$ of any point
$\lambda_x \in \g g_{\mathbb C}^*$.

We will also use a $U_{\BB R}$-invariant 2-form $\tau_{\lambda}$ on
${\cal B}$ defined by the formula
$$
\tau_{\lambda}(X_x,Y_x) = \lambda_x([X,Y]),
$$
where $X_x, Y_x \in T_x{\cal B}$ are the tangent vectors at $x \in {\cal B}$
induced by $X, Y \in \g u_{\BB R}$ via differentiation of the
$U_{\BB R}$-action.
It is an important property of $\lambda_x$ and $\tau_{\lambda}$
that together they make a (non-homogeneous) differential form on ${\cal B}$
$$
\langle X, \lambda_x \rangle - \tau_{\lambda}, \qquad X \in \g g_{\BB C},
$$
which depends holomorphically on $X$ (in fact linearly) and which is
{\em equivariantly closed} with respect to $U_{\BB R}$
(see \cite{BGV}, \cite{GS} or \cite{L3} for the definition of equivariantly
closed forms). It follows that the forms
$$
e^{\langle X, \lambda_x \rangle - \tau_{\lambda}}
\qquad \text{and} \qquad
e^{-\langle X, \lambda_x \rangle + \tau_{\lambda}}
$$
are $U_{\BB R}$-equivariantly closed forms too.
For a differential form $\omega$ which is possibly non-homogeneous,
we denote by $\omega_{[k]}$ its homogeneous component of degree $k$.
Then the integral character formula (Theorem 3.8 in \cite{SchV2})
says that the value of the character
$\tilde\theta$ at $\phi \in \Omega_c^{top}(\g g_{\BB R})$ is given by
\begin{multline}  \label{intcharformula}
\theta (\phi) =
\frac 1{(2\pi i)^dd!} \int_{Ch({\cal F})}
(\hat\phi \circ \mu_{\lambda}) \cdot (-\sigma_{\cal B}+\tau_{\lambda})^d  \\
= \frac 1{(-2\pi i)^d} \int_{Ch({\cal F})} \Bigl( \int_{\g g_{\BB R}}
e^{\langle X, \mu_{\cal B}(\zeta) + \lambda_x \rangle
+ \sigma_{\cal B} - \tau_{\lambda}}
\wedge \phi(X) \Bigr)_{[2d]},  \\
\qquad X \in \g g_{\BB R},\: \zeta \in T^*{\cal B},
\end{multline}
where $\mu_{\lambda} = \mu_{\cal B} + \lambda_x$ and
$d = \dim_{\mathbb C} {\cal B}$.
This integral converges because the expression $\hat\phi \circ \mu_{\lambda}$
involves the Fourier transform of $\phi$
and decays rapidly along the support of $Ch({\cal F})$.

\separate

\begin{section}
{Setup}
\end{section}

Recall that $G_{\BB C}$ is a connected complex algebraic reductive group
which is defined over $\BB R$, $G_{\BB C}$ acts algebraically on a
smooth complex projective variety $M$.
And we are primarily interested in representations of a real reductive
group $G_{\BB R} \subset G_{\BB C}$ lying between the group of real points
$G_{\BB C}(\BB R)$ and the identity component $G_{\BB C}(\BB R)^0$.

We will be using the concepts of ${\cal D}$-modules and derived categories;
\cite{Bo} and \cite{KaScha} are good introductions to these subjects.
Let ${\cal O}_M$ denote the structure sheaf on $M$ and
let ${\cal D}_M$ denote the sheaf of linear differential operators on $M$
with algebraic coefficients; ${\cal D}_M$ acts on ${\cal O}_M$.
The definition of a quasi-equivariant ${\cal D}$-module can be found
in many different sources; for convenience, we copy the definition given in
\cite{KaSchm}.

Let $G_{\BB C}$ be a connected complex algebraic reductive group
acting algebraically on a smooth complex algebraic variety $M$,
$\g g_{\mathbb C} = Lie(G_{\BB C})$.
We write $\mu$ for the action morphism and $\pi$ for the projection map
$G_{\BB C} \times M \to M$:
$$
\mu, \pi: G_{\BB C} \times M \to M,
\quad \mu(g, x) = gx, \quad \pi(g, x) = x
\qquad g \in G_{\BB C},\: x \in M.
$$
We also consider three maps
$$
q_j: G_{\BB C} \times G_{\BB C} \times M \to G_{\BB C} \times M,
\qquad j=1,2,3,
$$
$$
q_1(g_1,g_2,x)=(g_1,g_2x), \qquad
q_2(g_1,g_2,x)=(g_1g_2,x), \qquad
q_3(g_1,g_2,x)=(g_2,x).
$$
Then we have the following identities:
\begin{align*}
\mu \circ q_1 = \mu \circ q_2 &: \quad (g_1,g_2,x) \mapsto g_1g_2x,  \\
\pi \circ q_2 = \pi \circ q_3 &: \quad (g_1,g_2,x) \mapsto x, \\
\mu \circ q_3 = \pi \circ q_1 &: \quad (g_1,g_2,x) \mapsto g_2x.
\end{align*}

\begin{df}
We denote by ${\cal O}_{G_{\mathbb C}} \boxtimes {\cal D}_M$ the subalgebra
${\cal O}_{G_{\mathbb C} \times M} \otimes_{\pi^{-1} {\cal O}_M}
\pi^{-1}{\cal D}_M$ of ${\cal D}_{G_{\mathbb C} \times M}$.
A {\em quasi-$G_{\BB C}$-equivariant ${\cal D}_M$-module} is a
${\cal D}_M$-module $\g M$ equipped with an
${\cal O}_{G_{\mathbb C}} \boxtimes {\cal D}_M$-linear
isomorphism $\beta: \mu^* \g M \,\tilde\longrightarrow\, \pi^* \g M$, such that
the composition of isomorphisms
$$
\begin{CD}
q_2^* \mu^* \g M  \simeq q_1^* \mu^* \g M @>{q_1^*\beta}>>
q_1^* \pi^* \g M  \simeq  q_3^* \mu^* \g M  @>{q_3^*\beta}>> q_3^* \pi^* \g M
\simeq q_2^* \pi^* \g M
\end{CD}
$$
coincides with
$$
\begin{CD}
q_2^* \mu^* \g M   @>{q_2^*\beta}>> q_2^* \pi^* \g M.
\end{CD}
$$
If $\beta$ is linear even over ${\cal D}_{G_{\mathbb C} \times M}$,
this reduces to the usual definition of a {\em $G_{\mathbb C}$-equivariant
${\cal D}_M$-module}.
\end{df}

Informally speaking, this definition can be interpreted as follows.
For each $g \in G_{\mathbb C}$, denote by $\mu_g :M \to M$ the translation
by $g$, i.e. $\mu_g: x \mapsto gx$.
Then $\beta$ consists of a family of isomorphisms of ${\cal D}_M$-modules
$\beta_g: \mu_g^* \g M \,\tilde\longrightarrow\, \g M$,
depending algebraically on $g$ and multiplicative in the variable $g$.

\begin{ex}  \label{line_bundle}
Let $({\bf E}, \nabla_{\bf E})$ be a $G_{\BB C}$-equivariant algebraic line
bundle over a $G_{\BB C}$-invariant open algebraic subset
$U \subset M$ with a $G_{\BB C}$-invariant algebraic flat connection
$\nabla_{\bf E}$.
Take $\g M$ to be the direct image of the sheaf of sections
${\cal O}({\bf E})$, under the inclusion map $U \hookrightarrow M$.
Then ${\cal D}_M$ acts on ${\cal O}({\bf E})$ via the flat
connection making ${\cal O}({\bf E})$ a quasi-$G_{\BB C}$-equivariant
${\cal D}_M$-module.

Each element $X \in \g g_{\mathbb C}$ acts on $\g M$ in two ways.
One way is by inducing the vector field $-X_M$ given by (\ref{vfield})
which in turn operates on $\g M$ via the connection.
The other way is by infinitesimal translation of the sections of the
$G_{\BB C}$-equivariant line bundle ${\bf E}$.
When these two actions of $\g g_{\mathbb C}$ coincide, the ${\cal D}_M$-module
$\g M$ is $G_{\mathbb C}$-equivariant.
\end{ex}

Let $\g M$ be a coherent, quasi-$G_{\BB C}$-equivariant ${\cal D}_M$-module.
Recall that, for $X \in \g g_{\BB C}$, $X_M$ denotes the vector field on $M$
given by (\ref{vfield}).
Following \cite{KaSchm} we get two different actions of $\g g_{\BB C}$ on
$\g M$. The first action is via
$\g g_{\BB C} \ni X \mapsto -X_M \in \Gamma({\cal D}_M)$
-- the global sections of ${\cal D}_M$ -- followed by the ${\cal D}_M$-module
structure; we denote this action by
$\alpha_{\cal D}$.
And the second action is through differentiation of the $G_{\BB C}$-action
when we regard $\g M$ as a $G_{\BB C}$-equivariant ${\cal O}_M$-module;
this action is denoted by $\alpha_t$.
We set $\gamma = \alpha_t - \alpha_{\cal D}$, then $\gamma$ is a Lie algebra
homomorphism
\begin{equation}  \label{gamma}
\gamma: \g g_{\BB C} \longrightarrow \End_{{\cal D}_M} (\g M).
\end{equation}
This way $\g M$ becomes a $({\cal D}_M, {\cal U}(\g g_{\BB C}))$-module,
where ${\cal U}(\g g_{\BB C})$ denotes the universal enveloping algebra
of $\g g_{\BB C}$.
The quasi-$G_{\BB C}$-equivariant ${\cal D}_M$-module $\g M$ is
$G_{\BB C}$-equivariant precisely when $\gamma \equiv 0$.

We say that the quasi-$G_{\BB C}$-equivariant ${\cal D}_M$-module $\g M$ is
{\em ${\cal Z}(\g g_{\BB C})$-finite} if some ideal of finite codimension
${\cal I} \subset {\cal Z}(\g g_{\BB C})$ (the center of 
${\cal U}(\g g_{\BB C})$) annihilates $\g M$ via the $\gamma$-action.

We denote by $Ch(\g M)$ the characteristic cycle of
$\g M$ which lies in $T^*M$.
Pick a Borel subalgebra $\g b_{\BB C} \subset \g g_{\BB C}$
and define a subset of $\g g_{\BB C}^*$
$$
\g b_{\BB C}^{\perp} =
\{ \xi \in \g g_{\BB C}^* ;\: \xi|_{\g b_{\BB C}} \equiv 0 \}.
$$
Then the ${\cal D}_M$-module $\g M$ is called {\em admissible} if
$$
Ch(\g M) \cap \mu_M^{-1} (\g b_{\BB C}^{\perp}) \quad \subset \quad T^*M
$$
is a Lagrangian variety.
When $\g M$ is ${\cal Z}(\g g_{\BB C})$-finite,
the variety $Ch(\g M) \cap \mu_M^{-1} (\g b_{\BB C}^{\perp})$
is known to be involutive (\cite{KaMF} or \cite{Gi}),
hence the above condition is equivalent to
$$
\dim_{\BB C} \bigl( Ch(\g M) \cap \mu_M^{-1} (\g b_{\BB C}^{\perp}) \bigr)
= \dim_{\BB C} M.
$$
Because of $G_{\BB C}$-invariance of $Ch(\g M)$, if this condition is
satisfied for one particular Borel subalgebra
$\g b_{\BB C} \subset \g g_{\BB C}$,
then it is satisfied for all Borel subalgebras of $ \g g_{\BB C}$,
and this definition is independent of the choice of the Borel subalgebra
$\g b_{\BB C}$.

\separate

Let ${\cal S}$ be a $G_{\BB R}$-equivariant constructible sheaf on $M$.
We denote by ${\cal O}_{M^{an}}$ the sheaf of holomorphic functions on $M$.
Then M.~Kashiwara and W.~Schmid \cite{KaSchm} equip the vector spaces
\begin{equation} \label{ext}
\rhom^p_{{\cal D}_M} ( \g M \otimes {\cal S}, {\cal O}_{M^{an}}),
\qquad p \in \BB Z,
\end{equation}
with a natural Fr\'echet topology and prove that the resulting virtual
representation of $G_{\BB R}$
\begin{equation}  \label{vrep}
\sum_p (-1)^p
\rhom^p_{{\cal D}_M} ( \g M \otimes {\cal S}, {\cal O}_{M^{an}})
\end{equation}
is admissible of finite length whenever the ${\cal D}_M$-module $\g M$
is admissible and ${\cal Z}(\g g_{\BB C})$-finite.
In particular, the representation (\ref{vrep}) has a character $\theta$
in the sense of Harish-Chandra. As before, $\theta$ is a distribution on
$\g g_{\mathbb R}$.

In this paragraph we outline M.~Kashiwara and W.~Schmid's construction of
topology on the spaces (\ref{ext}). Oversimplifying and ignoring the
$G_{\mathbb C}$- and $G_{\mathbb R}$-actions, suppose first that the sheaf
${\cal S}$ is $(j_{U \hookrightarrow M})_! \mathbb C_U$,
where $U \subset M$ is an open semi-algebraic $G_{\mathbb R}$-invariant subset
and $\mathbb C_U$ is the constant sheaf on $U$.
Furthermore, suppose that $\g M$ is a {\em locally free}
quasi-$G_{\mathbb C}$-equivariant ${\cal D}_M$-module, i.e.
$\g M = {\cal D}_M \otimes_{{\cal O}_M} \g F$ for some coherent,
locally free, $G_{\mathbb C}$-equivariant ${\cal O_M}$-module $\g F$ on $M$.
They replace the sheaf ${\cal O}_{M^{an}}$ with the ${\cal C}^{\infty}$
Dolbeault complex $\Omega_M^{(0,\cdot)}$, to which it is quasi-isomorphic,
and write out isomorphisms of complexes of vector spaces without any topology:
\begin{multline}  \label{topology}
\rhom_{{\cal D}_M} ( \g M \otimes {\cal S}, {\cal O}_{M^{an}})
\simeq \rhom_{{\cal D}_M} ( {\cal D}_M \otimes_{{\cal O}_M} \g F
\otimes (j_{U \hookrightarrow M})_! \mathbb C_U, \Omega_M^{(0,\cdot)})  \\
\simeq \rhom_{{\cal O}_M} ( \g F \otimes
(j_{U \hookrightarrow M})_! \mathbb C_U, \Omega_M^{(0,\cdot)})
\simeq \operatorname{R}\Gamma ( U;
(\g F^*)^{an} \otimes_{{\cal O}_{M^{an}}} \Omega_M^{(0,\cdot)})  \\
\simeq \Gamma ( U;
(\g F^*)^{an} \otimes_{{\cal O}_{M^{an}}} \Omega_M^{(0,\cdot)}),
\end{multline}
where $(\g F^*)^{an}$ denotes the sheaf of analytic sections of the dual
of $\g F$. The complex on the right has a natural Fr\'echet topology -- the
${\cal C}^{\infty}$ topology for differential forms -- and continuous
$G_{\mathbb R}$ action.
General $\g M$ and ${\cal S}$, have resolutions by locally free
quasi-$G_{\mathbb C}$-equivariant ${\cal D}_M$-modules
(Lemma 4.7 in \cite{KaSchm}) and by sheaves of the type
$(j_{U \hookrightarrow M})_! \mathbb C_U$ respectively.
This, combined with the acyclicity of the right hand side of (\ref{topology}),
makes it possible to equip the vector spaces (\ref{ext}) with a
Fr\'echet topology. Then M.~Kashiwara and W.~Schmid work hard to show that
this topology does not depend on the choices involved. They do it by
proving that the topology on the spaces (\ref{ext}) is that of the
maximal globalization of their underlying Harish-Chandra modules.

\separate

Recall that ${\cal B}$ denotes the flag variety of $\g g_{\BB C}$.
We will establish a formula for this character as an integral over a
cycle in $T^*({\cal B} \times M)$ under an additional assumption that
$\g M$ has an infinitesimal character, i.e. ${\cal Z}(\g g_{\BB C})$
acts on $\g M$ by a character.
Even if this condition is not satisfied, since $\g M$ is assumed to be 
${\cal Z}(\g g_{\BB C})$-finite, there is a finite filtration of $\g M$
by quasi-$G_{\BB C}$-equivariant ${\cal D}_M$-submodules such that the
successive quotients have infinitesimal characters, and we can apply this
integral character formula to each of these quotients separately.

In the special case when $G_{\BB C} = \BB C^{\times}$, the
${\cal D}_M$-module $\g M$ is the sheaf of sections of a
$G_{\BB C}$-equivariant line bundle $({\bf E}, \nabla_{\bf E})$ over $M$
as in Example \ref{line_bundle} with $U=M$,
and ${\cal S}$ is the constant sheaf $\BB C_M$,
the flag variety ${\cal B}$ consists of just one  point and the integral
character formula will reduce to the classical Riemann-Roch-Hirzebruch formula.

\begin{rem}
We do not need $M$ to be projective to establish that (\ref{vrep}) is
admissible of finite length. In fact it is sufficient to assume that $M$ is
a smooth {\em quasi-projective} variety.
The compactness of $M$ will be needed for the integral character formula.

On the other hand, the result of Sumihiro \cite{Su} restated as
Proposition 4.6 in \cite{KaSchm} together with Theorem 5.12 there
show that we can always embed $M$ into a smooth projective variety
and there is no loss of generality in assuming that $M$ is projective.
We will illustrate this in Section \ref{example}.
\end{rem}

Our derivation of the integral character formula for $\theta$ follows the
following scheme:
\begin{itemize}
\item
Replace the pair $(\g M, {\cal S})$ which lives on $M$ with a new pair
$(\widetilde{\g M}, \widetilde{\cal S})$ which lives on the flag variety
${\cal B}$ such that the virtual representation (\ref{vrep})
stays unchanged (\ref{newpair});
\item
Write out the fixed point character formula for $\theta$ in terms of
$(\widetilde{\g M}, \widetilde{\cal S})$ and zeroes on ${\cal B}$;
\item
We want to prove that the integral (\ref{intcharformula}) represents $\theta$;
first we apply the integral localization formula \cite{L3} to
(\ref{intcharformula});
\item
Combine similar terms in the result of the previous step and match them with
the terms in the fixed point character formula for $\theta$ obtained in the
earlier step.
\end{itemize}

\separate

\begin{section}
{Statement of the main result}
\end{section}

In the previous section we assumed that $\g M$ was a 
quasi-$G_{\BB C}$-equivariant coherent ${\cal D}_M$-module on $M$
which was admissible and had an infinitesimal character.
We will follow \cite{KaSchm} and index characters of ${\cal Z}(\g g_{\BB C})$
by linear functionals $\lambda$ on the universal Cartan without the
customary shift by $\rho$ (half sum of the positive roots);
in other words, $\chi_{\lambda}: {\cal Z}(\g g_{\BB C}) \to \BB C$
denotes the character by which ${\cal Z}(\g g_{\BB C})$ acts on the Verma
module with highest weight $\lambda$. Then $\chi_{\lambda} = \chi_{\mu}$
if and only if $\lambda + \rho$ is conjugate to $\mu + \rho$ under the
action of the Weyl group $W$.

\begin{rem}  \label{twistings}
In this article we use results from \cite{KaSchm} and \cite{SchV2}.
Unfortunately these two sources use different conventions for labeling
characters and twists. We follow the conventions of \cite{KaSchm} explained
above. On the other hand \cite{SchV2} shift these notations by $\rho$ and
they define, for instance, the ``$G_{\BB R}$-equivariant derived category on
${\cal B}$ with twist $(\lambda-\rho)$'' denoted by
$\operatorname{D}_{G_{\BB R}} ({\cal B})_{\lambda}$
which becomes the bounded $G_{\BB R}$-equivariant derived category in the
usual ``untwisted'' sense precisely when $\lambda = \rho$.
Similarly, ${\cal O}_{\cal B}(\lambda)$ in \cite{SchV2} denotes
the twisted sheaf of holomorphic functions with twist $(\lambda-\rho)$
and which becomes the sheaf of ordinary holomorphic functions precisely
when $\lambda = \rho$.
\end{rem}

We assume that $\g M$ is an object in $\modg^{coh,\lambda-\rho}({\cal D}_M)$
-- the category of quasi-$G_{\BB C}$-equivariant coherent
${\cal D}_M$-modules on $M$ with infinitesimal character $\chi_{\lambda-\rho}$.

The Kashiwara-Schmid construction of Fr\'echet topology on the spaces
(\ref{ext}) was carried out on the level of derived categories.
This means that instead of a sheaf ${\cal S}$ we have an element
${\cal S} \in D^b_{G_{\BB R}, \BB R-c} (\BB C_M)$ -- the bounded
derived category of $G_{\BB R}$-equivariant $\BB R$-constructible
sheaves, and we have a pairing
$$
\modg^{coh,\lambda-\rho}({\cal D}_M) \times D^b_{G_{\BB R}, \BB R-c} (\BB C_M)
\quad \longrightarrow \quad D^b({\cal F}_{G_{\BB R}}),
$$
$$
(\g M, {\cal S}) \quad \mapsto \quad
\rhom^{top}_{{\cal D}_M} ( \g M \otimes {\cal S}, {\cal O}_{M^{an}}),
$$
where $D^b({\cal F}_{G_{\BB R}})$ denotes the derived category
of $G_{\BB R}$-representations defined in \parag 3 of \cite{KaSchm}.
The category $D^b({\cal F}_{G_{\BB R}})$ is built on complexes of topological
vector spaces with continuous $G_{\mathbb R}$-actions $(C,d_C)$,
the differential maps $d_C^n : C^n \to C^{n+1}$ are required to be continuous
and $G_{\mathbb R}$-equivariant. A complex $(C,d_C)$ is {\em exact} if,
for all $n$,
$$
d_C^n : C^n \to \ker d_C^{n+1}
$$
is onto and is an open map relative to the subspace topology on
$\ker d_C^{n+1}$.
A complex $(C,d_C)$ becomes zero in $D^b({\cal F}_{G_{\BB R}})$ precisely
when it is exact.

Because we will apply results of A.~Beilinson and J.~Bernstein on equivalences
of categories \cite{BB}, we will also assume that $\lambda$ is
{\em integrally dominant}, i.e.
$$
\langle \check \alpha, \lambda \rangle \notin \BB Z_{<0},
\qquad \text{for every positive coroot $\check \alpha$.}
$$
This can always be achieved by replacing $\lambda$ with an appropriate
$W$-translate.

\separate

We denote by $\modge({\cal D}_{{\cal B},\,\lambda-\rho})$ the category of
modules over the sheaf of twisted differential operators
${\cal D}_{{\cal B},\,\lambda-\rho}$ on the flag variety ${\cal B}$.
The sheaf ${\cal D}_{{\cal B},\,\lambda-\rho}$ is defined in \cite{BB},
but we follow the twisting conventions of \cite{KaSchm}
as explained in Remark \ref{twistings}, so that
${\cal D}_{{\cal B},\,0} = {\cal D}_{{\cal B}}$ -- the sheaf of differential
operators on ${\cal B}$ without any twist.
We form a product space ${\cal B} \times M$ with diagonal $G_{\BB C}$-action
and consider the sheaf of twisted differential operators
${\cal D}_{{\cal B} \times M,\,\lambda-\rho} =_{def}
{\cal D}_{{\cal B},\,\lambda-\rho} \boxtimes {\cal D}_M$,
the twisting is confined to the factor ${\cal B}$.
We also denote by $\modge({\cal D}_{{\cal B} \times M,\,\lambda-\rho})$
the category of $G_{\BB C}$-equivariant
${\cal D}_{{\cal B} \times M,\,\lambda-\rho}$-modules on ${\cal B} \times M$.
Let $p$ and $q$ be the projection maps
$$
\begin{CD}
{\cal B} @<p<< {\cal B} \times M \\ & & @VVqV \\ & & M,
\end{CD}
$$
and let $\tilde p$ and $\tilde q$ be the induced projections on the
cotangent bundles:
$$
\begin{CD}
T^*{\cal B} @<\tilde p<< T^*({\cal B} \times M) \\ & & @VV\tilde qV \\
& & T^*M.
\end{CD}
$$

M.~Kashiwara and W.~Schmid \cite{KaSchm} use results of
A.~Beilinson and J.~Bernstein on equivalence of categories \cite{BB}
to prove that the pair $(\g M, {\cal S})$ on $M$ can be replaced with a
${\cal D}_{{\cal B},\,\rho-\lambda}$-module and a complex of sheaves on
${\cal B}$.
We will make this statement precise.
First of all, they show that there exists a coherent holonomic module
$\g L \in \modge({\cal D}_{{\cal B} \times M,\,\lambda-\rho})$ such that
$$
q_* \g L = \g M  \qquad \text{and} \qquad
R^k q_* \g L = 0 \text{ if $k \ne 0$,}
$$
namely one can take
$$
\g L = {\cal D}_{{\cal B},\,\lambda-\rho} \boxtimes \g M /
\gamma_{{\cal B} \times M}
(\g g_{\BB C})({\cal D}_{{\cal B},\,\lambda-\rho} \boxtimes \g M),
$$
where
$\gamma_{{\cal B} \times M}
(\g g_{\BB C})({\cal D}_{{\cal B},\,\lambda-\rho} \boxtimes \g M)$
denotes the image in ${\cal D}_{{\cal B},\,\lambda-\rho} \boxtimes \g M$ of
the module
$\g g_{\BB C} \otimes ({\cal D}_{{\cal B},\,\lambda-\rho} \boxtimes \g M)$
under the map $\gamma_{{\cal B} \times M}$ given by the equation (\ref{gamma})
with the ambient manifold ${\cal B} \times M$.

Then we apply the twisted deRham functor to obtain a complex of sheaves
$$
{\cal L} = DR_{{\cal B} \times M} (\g L) =_{def}
R{\cal H}om_{{\cal D}_{{\cal B} \times M, \, \lambda-\rho}}
({\cal O}_{\cal B}(\lambda-\rho) \boxtimes {\cal O}_M, \g L),
$$
which is an element of the bounded $G_{\BB C}$-equivariant,
$\BB C$-constructible derived category with twist $(\rho-\lambda)$ along the
${\cal B}$-factor denoted
$D^b_{G_{\BB C},\, \rho-\lambda, \BB C-c} (\BB C_{{\cal B} \times M})$.
Here ${\cal O}_{\cal B}(\lambda-\rho)$ is a twisted sheaf of holomorphic
functions on ${\cal B}$, with twist $(\lambda-\rho)$
(so that ${\cal O}_{\cal B}(0)$ is just the sheaf of functions on ${\cal B}$
with no twist at all).

Let $d = \dim_{\BB C} {\cal B}$, and let
$(\Omega_{\cal B}^d)^{-1} \simeq {\cal O}_{\cal B} (2\rho)$ denote the
reciprocal of the canonical sheaf.
Combining the equations (6.7), (7.14) and (7.15) from \cite{KaSchm}
we obtain
\begin{multline}  \label{maineqn}
\rhom^{top}_{{\cal D}_M} ( \g M \otimes {\cal S}, {\cal O}_{M^{an}})  \\
\simeq
\rhom^{top}_{{\cal D}_{\cal B},\, \rho-\lambda}
\bigl( {\cal D}_{{\cal B},\, \rho-\lambda} \otimes_{{\cal O}_{\cal B}}
(\Omega_{\cal B}^d)^{-1} \otimes Rp_*({\cal L} \otimes q^{-1}{\cal S}),
{\cal O}_{{\cal B}^{an}}(\rho-\lambda) \bigr) [d - 2\dim_{\BB C} M]  \\
\simeq
\rhom^{top}_{{\cal D}_{\cal B}, -\lambda -\rho}
\bigl( {\cal D}_{{\cal B}, -\lambda -\rho}
\otimes Rp_*({\cal L} \otimes q^{-1}{\cal S}),
{\cal O}_{{\cal B}^{an}}(-\lambda -\rho) \bigr) [d - 2\dim_{\BB C} M]  \\
\simeq
\rhom \bigl( Rp_*({\cal L} \otimes q^{-1}{\cal S}), 
{\cal O}_{{\cal B}^{an}}(-\lambda -\rho) \bigr) [d - 2\dim_{\BB C} M]
\end{multline}
as elements of the derived category of $G_{\BB R}$-representations
$D^b({\cal F}_{G_{\BB R}})$.
Here we view $Rp_*({\cal L} \otimes q^{-1}{\cal S})$ as an object in
$D^b_{G_{\BB R},-\lambda-\rho, \BB R-c} (\BB C_{\cal B})$,
which makes sense because $2\rho$ is an integral weight, and this implies
the existence of a canonical isomorphism (6.8) in \cite{KaSchm}
$$
D^b_{G_{\BB R},\, \rho-\lambda, \BB R-c} (\BB C_{\cal B})
\quad\simeq\quad
D^b_{G_{\BB R},-\lambda-\rho, \BB R-c} (\BB C_{\cal B}).
$$
That is, the pair $(\g M, {\cal S})$ on the variety $M$ is replaced by a pair
\begin{equation}  \label{newpair}
\bigl( {\cal D}_{{\cal B}, -\lambda -\rho},
Rp_*({\cal L} \otimes q^{-1}{\cal S}) \bigr)
\end{equation}
on the flag variety ${\cal B}$, with an additional twist by $(-\lambda -\rho)$.

\separate

Our starting point is the integral character formula (\ref{intcharformula})
proved by W.~Schmid and K.~Vilonen in \cite{SchV2}. We apply it to the
right hand side of (\ref{maineqn}).
We fix a compact real form $U_{\BB R} \subset G_{\BB C}$.
Recall that $\mu_{\cal B}: T^*{\cal B} \to \g g_{\BB C}^*$ is the moment
map defined by (\ref{mu}), $\sigma_{\cal B}$ is the canonical complex
algebraic holomorphic symplectic form on $T^*{\cal B}$,
$\lambda_x : {\cal B} \to \g g_{\BB C}^*$ is the $U_{\mathbb R}$-equivariant
Rossmann's map, and $\tau_{\lambda}$ is a certain $U_{\BB R}$-invariant
2-form on ${\cal B}$.
The character of the virtual representation (\ref{vrep}) is a distribution
on $\Omega_c^{top}(\g g_{\BB R})$ -- the space of smooth compactly
supported differential forms on $\g g_{\BB R}$ of top degree.
For an element $\phi \in \Omega_c^{top}(\g g_{\BB R})$, its Fourier transform 
$\hat \phi \in {\cal C}^{\infty}(\g g_{\BB C}^*)$ is defined by (\ref{FT})
without the customary factor of $i = \sqrt{-1}$ in the exponent.
Then the integral character formula says that the character
$\theta$ of the virtual representation (\ref{maineqn}) of $G_{\BB R}$,
as a distribution on $\Omega_c^{top}(\g g_{\BB R})$, is
\begin{multline}  \label{intcharformula1}
\theta(\phi) = \frac 1{(2\pi i)^d d!}
\int_{Ch(Rp_*({\cal L} \otimes q^{-1}{\cal S}))^a}
(\hat\phi \circ \mu_{-\lambda}) \cdot (\sigma_{\cal B} + \tau_{\lambda})^d  \\
= \frac 1{(2\pi i)^d}
\int_{Ch(Rp_*({\cal L} \otimes q^{-1}{\cal S}))^a}
\Bigl( \int_{\g g_{\BB R}}
e^{\langle X, \mu_{\cal B}(\zeta) - \lambda_x \rangle
+ \sigma_{\cal B} + \tau_{\lambda}}
\wedge \phi(X) \Bigr)_{[2d]},  \\
\qquad X \in \g g_{\BB R},\: \zeta \in T^*{\cal B}.
\end{multline}
Here $\mu_{-\lambda} = \mu_{\cal B} - \lambda_x$,
$d= \dim_{\BB C} {\cal B}$, $a: T^*{\cal B} \to T^*{\cal B}$
is the antipodal map $\zeta \mapsto -\zeta$, and
$Ch(Rp_*({\cal L} \otimes q^{-1}{\cal S}))^a$ denotes
the image under this antipodal map of the characteristic cycle of 
$Rp_*({\cal L} \otimes q^{-1}{\cal S})$ (which is a cycle in $T^*{\cal B}$).

\begin{rem}
If $Z$ is a complex manifold and
$Z^{\BB R}$ is the underlying real analytic manifold,
then the holomorphic symplectic form $\sigma_Z$ is defined on the
holomorphic cotangent bundle $T^*Z$, while the characteristic cycles
of constructible sheaves on $Z$ lie in the real cotangent bundle
$T^*(Z^{\BB R})$. Hence we need to identify $T^*Z$ with $T^*(Z^{\BB R})$.
There are at least two different but equally natural ways of doing this,
we use the convention (11.1.2) of \cite{KaScha}, Chapter XI;
the same convention is used in \cite{L1}, \cite{L2}, \cite{L3} and
\cite{SchV2}. Under this convention, if
$\sigma_{Z,\BB R}$ is the canonical real symplectic form on $T^*Z^{\BB R}$
and $\sigma_Z$ is the canonical complex holomorphic symplectic form on $T^*Z$,
then $\sigma_{Z,\BB R}$ gets identified with $2 \re \sigma_Z$.
\end{rem}

Another important ingredient is a generalization of the Hopf index theorem
stated as Corollary 9.5.2 in \cite{KaScha}.
Let ${\cal T}$ be a constructible sheaf on $M$ or
an element of the derived category $D^b_{\BB R-c} (\BB C_M)$,
and let $\chi(M, {\cal T})$ denote the Euler characteristic of $M$ with
respect to ${\cal T}$.
Then
\begin{equation}  \label{hopf}
\chi(M, {\cal T}) = \# \bigl( [M] \cap Ch({\cal T}) \bigr),
\end{equation}
where $[M]$ denotes the fundamental cycle of $M$.
Alternatively one can apply the equation (5.30) from \cite{Schu}.
Let $\thom_{T^*M}$ denote the Thom form of the cotangent bundle
$T^*M \twoheadrightarrow M$.
That is $\thom_{T^*M}$ is a closed differential form on $T^*M$ of degree
$2\dim_{\BB C}M$ which decays rapidly along the fiber (or even compactly
supported along the fiber) and such that
$$
\int_{T_x^*M} \thom_{T^*M} = (2\pi)^{\dim_{\BB C}M}, \qquad \forall x \in M.
$$
We regard $M$ as a submanifold of $T^*M$ via the zero section inclusion.
Then the restriction of the Thom form to $M$ is the Euler form of $M$.
Since the form $(2\pi)^{-\dim_{\BB C}M} \thom_{T^*M}$ is Poincar\'e dual to
the homology class of $[M]$ in $T^*M$, the Hopf index theorem (\ref{hopf})
can be rewritten as
\begin{equation}  \label{bonnet}
\chi(M, {\cal T}) = \# \bigl( [M] \cap Ch({\cal T}) \bigr)
= (2\pi)^{-\dim_{\BB C}M} \int_{Ch({\cal T})} \thom_{T^*M},
\end{equation}
which is a generalization of the Gauss-Bonnet theorem.

Recall that $U_{\BB R}$ is a compact real form of $G_{\BB C}$.
The form $\thom_{T^*M}$ may be chosen to be $U_{\BB R}$-invariant.
If, in addition, the cotangent bundle $T^*M$ has a spin structure,
then V.~Mathai and D.~Quillen showed in \cite{MQ}
(see also Section 7.7 in \cite{BGV})
that $\thom_{T^*M}$ can be realized as the top degree part of a
$U_{\BB R}$-equivariantly closed form on $T^*M$ in a canonical way.

\separate

The last essential ingredient is the deformation argument for integrals
of equivariantly closed forms from \cite{L3}.
This argument requires that any maximal complex torus
$T_{\BB C} \subset G_{\BB C}$ acts on $M$ with isolated fixed points.
(Since $M$ is compact, there will be only finitely many of those.)
If this condition is satisfied for one particular torus, then it is
satisfied for all tori because all maximal complex tori are conjugate
by elements of $G_{\BB C}$.
(This condition is satisfied in all application we have in mind.)

We will combine the integral character formula (\ref{intcharformula1})
and the Gauss-Bonnet formula (\ref{bonnet}) to get a formula for the
character of the virtual representation (\ref{maineqn}) as an integral over
a cycle in $T^*({\cal B} \times M)$.
Set
$$
\widetilde{\cal S} = {\cal L} \otimes q^{-1}{\cal S}
\qquad \text{and} \qquad
\Lambda = Ch(\widetilde{\cal S})^a
= Ch(\BB D_{{\cal B} \times M} (\widetilde{\cal S}))
\subset T^*({\cal B} \times M).
$$
Here $a: T^*({\cal B} \times M) \to T^*({\cal B} \times M)$
is the antipodal map $\zeta \mapsto -\zeta$ and
$Ch(\widetilde{\cal S})^a$ denotes the image under this antipodal map of
the characteristic cycle of $\widetilde{\cal S}$; the operator
$\BB D_{{\cal B} \times M}$ is the Verdier duality operator
(see (\ref{D_XxM}) below for its properties).

\begin{thm}  \label{main}
Suppose that any maximal complex torus $T_{\BB C} \subset G_{\BB C}$ acts on
$M$ with isolated fixed points.
Then the value of the character $\theta$ of the virtual representation
(\ref{vrep}) on $\phi \in \Omega_c^{top}(\g g_{\BB R})$ is
\begin{multline}  \label{mainintegral}
\theta(\phi) =
\frac {i^n}{(2\pi i)^{d+n} d!}
\int_{\Lambda}
\tilde p^* \bigl( (\hat\phi \circ \mu_{-\lambda}) \cdot
(\sigma_{\cal B} + \tau_{\lambda})^d \bigr)
\wedge \tilde q^*\thom_{T^*M}  \\
= \frac {i^n}{(2\pi i)^{d+n}}
\int_{\Lambda} \Bigl( \int_{\g g_{\BB R}} 
\tilde p^* e^{\langle X, \mu_{\cal B}(\zeta) - \lambda_x \rangle +
\sigma_{\cal B} + \tau_{\lambda}}
\wedge \tilde q^*\thom_{T^*M} \wedge \phi(X) \Bigr)_{[2(d+n)]},
\end{multline}
where 
$\mu_{-\lambda} = \mu_{\cal B} - \lambda_x$,
$d =\dim_{\BB C} {\cal B}$ and $n=\dim_{\BB C} M$.
\end{thm}

\begin{rem}
Suppose that the group $G_{\BB C} = \BB C^{\times}$, the
${\cal D}_M$-module $\g M$ is the sheaf of sections of a
$G_{\BB C}$-equivariant line bundle $({\bf E}, \nabla_{\bf E})$
over $M$ as in Example \ref{line_bundle},
and ${\cal S}$ is the constant sheaf $\BB C_M$.
Then the flag variety ${\cal B}$ is just one  point,
the cycle $\Lambda =[M]$, and because of the Riemann-Roch relationship
(8.4) in \cite{MQ} (or Theorem 7.44 in \cite{BGV})
the above integral character formula will reduce to the classical
Riemann-Roch-Hirzebruch formula.
There is no curvature of the line bundle $({\bf E}, \nabla_{\bf E})$
present in the character formula (\ref{mainintegral})
because we assume that the connection $\nabla_{\bf E}$ is flat so that
${\cal O}({\bf E})$ is a quasi-$G_{\BB C}$-equivariant ${\cal D}_M$-module.
\end{rem}

\separate

\begin{section}
{Proof of Theorem \ref{main}}
\end{section}

We start our proof by applying the fixed point character formula
(\ref{fpcharformula}) combined with (\ref{m})
due to W.~Schmid and K.~Vilonen \cite{SchV2} to the right hand side
of (\ref{maineqn}).
Thus the character $\theta$ of the virtual representation (\ref{maineqn})
is given by integration against a function $F_{\theta}$ on $\g g_{\BB R}$:
\begin{equation}  \label{F}
\theta(\phi) = \int_{\g g_{\BB R}} F_{\theta} \phi,
\qquad \phi \in \Omega_c^{top}(\g g_{\BB R}).
\end{equation}
This function $F_{\theta}$ is an $Ad(G_{\BB R})$-invariant, locally $L^1$
function on $\g g_{\BB R}$ whose restriction to the set of regular semisimple
elements $\g g_{\BB R}^{rs}$ can be represented by a real analytic function.
The value of this analytic function at $X \in \g g_{\BB R}^{rs}$ is determined
by zeroes of the vector field $X_{\cal B}$ on the flag variety ${\cal B}$
as follows.

Recall that $\g t_{\BB C}(X) \subset \g g_{\BB C}$ is the unique Cartan
subalgebra containing $X \in \g g_{\BB R}^{rs}$.
Let $\Psi(X) \subset \g t^*_{\BB C}(X)$ be the root system of $\g g_{\BB C}$
with respect to $\g t_{\BB C}(X)$ and pick a positive root system
$\Psi^{\le}(X) \subset \Psi$ such that
$$
\re \alpha (X) \le 0 \qquad \text{for all $\alpha \in \Psi^{\le}(X)$.}
$$
Let $B(X) \subset G_{\BB C}$ be the Borel subgroup whose Lie algebra
contains $\g t_{\BB C}(X)$ and the positive root spaces corresponding
to $\Psi^{\le}(X)$.
The action of $B(X)$ on ${\cal B}$ has exactly the same number of orbits
$O_1, \dots, O_{|W|}$ as the order of the Weyl group $W$ of $G_{\BB C}$.
Each of these orbits $O_k$ contains exactly one zero of the vector field
$X_{\cal B}$, and we order the zeroes $\{\g b_1,\dots,\g b_{|W|}\}$ of
$X_{\cal B}$ so that $\g b_k$ is contained in $O_k$, $k=1,\dots,|W|$.
The set of zeroes $\{\g b_1,\dots,\g b_{|W|}\}$ is precisely the set of Borel
subalgebras containing $\g t_{\BB C}(X)$.
Let $\Psi^+_{\g b_k}(X) \subset \Psi(X)$ be the
positive root system such that $\g b_k$ contains all the negative root spaces,
$k=1,\dots,|W|$.
Then the fixed point character formula \cite{SchV2} says that
the function $F_{\theta}$ which appeared in the equation (\ref{F}) is
\begin{equation}  \label{fpcharf}
F_{\theta} (X) = (-1)^d \sum_{k=1}^{|W|} m_{\g b_k}(X)
\frac {e^{-\langle X , \lambda_{\g b_k} \rangle}}
{\prod_{\alpha \in \Psi^+_{\g b_k}(X)} \alpha(X)},
\end{equation}
where $m_{\g b_k}(X)$'s are integer multiplicities given by the formula
\begin{equation}  \label{m1}
m_{\g b_k}(X) = \chi
\bigl( {\cal H}^*_{O_k} (Rp_*(\widetilde{\cal S}))_{\g b_k} \bigr)
= \chi \bigl( (j_{\{\g b_k\} \hookrightarrow O_k})^* \circ
(j_{O_k \hookrightarrow {\cal B}})^! \bigl(
Rp_*(\widetilde{\cal S}) \bigr) \bigr).
\end{equation}

Let $\BB D_{\cal B}$ and $\BB D_{{\cal B} \times M}$ denote the
Verdier duality operators on ${\cal B}$ and ${\cal B} \times M$ respectively:
$$
\BB D_{\cal B} : \quad
D^b_{G_{\BB R},-\lambda-\rho,\BB R-c} (\BB C_{\cal B})
\quad \tilde\longrightarrow \quad
D^b_{G_{\BB R},\lambda+\rho,\BB R-c} (\BB C_{\cal B});
$$
\begin{equation} \label{D_XxM}
\BB D_{{\cal B} \times M} : \quad
D^b_{G_{\BB R},-\lambda-\rho,\BB R-c} (\BB C_{{\cal B} \times M})
\quad \tilde\longrightarrow \quad
D^b_{G_{\BB R},\lambda+\rho,\BB R-c} (\BB C_{{\cal B} \times M}).
\end{equation}
The effect of the Verdier duality operator $\BB D_Z$ on the
characteristic cycle of an $\BB R$-constructible sheaf ${\cal T}$
(or an element of $D^b_{\BB R-c} (\BB C_Z)$) on any
smooth quasi-projective variety $Z$ is described by
$$
Ch(\BB D_Z ({\cal T})) = Ch({\cal T})^a,
$$
where $a: T^*Z \to T^*Z$ is the antipodal map $\zeta \mapsto -\zeta$ and
$Ch({\cal T})^a$ denotes the image under this antipodal map of
the characteristic cycle of ${\cal T}$.

\separate

For a regular semisimple element $X \in \g g_{\BB C}^{rs}$ we denote by
$T_{\BB C}(X) = \exp (\g t_{\BB C}(X)) \subset G_{\BB C}$
the maximal complex torus corresponding to the unique Cartan
subalgebra containing $X$.
If $x \in M$ is a point fixed by $T_{\BB C}(X)$; then $T_{\BB C}(X)$ acts
linearly on the tangent space $T_xM$.
We denote by $\Delta (X) \subset \g t_{\BB C}(X)^*$ the set of weights
of $\g t_{\BB C}(X)$ which are either roots of $\g g_{\BB C}$ or occur in
the tangent space $T_xM$ of some point $x \in M$ fixed by $T_{\BB C}(X)$.

Let $\g g_{\BB R}'$ denote the set of {\em strongly regular semisimple}
elements $X \in \g g_{\BB R}^{rs}$ which satisfy the following additional
properties. If $\g t_{\BB R}(X) \subset \g g_{\BB R}$ and
$\g t_{\BB C}(X) \subset \g g_{\BB C}$ are the unique Cartan
subalgebras in $\g g_{\BB R}$ and $\g g_{\BB C}$ respectively
containing $X$, then:
\begin{enumerate}
\item
The set of zeroes of the vector field $X_M$ is exactly the set of points in
$M$ fixed by the complex torus
$T_{\BB C}(X) = \exp (\g t_{\BB C}(X)) \subset G_{\BB C}$;
\item
$\beta(X) \ne 0$ for all $\beta \in \Delta(X) \subset \g t_{\BB C}(X)^*$;
\item
For each $\beta \in \Delta(X)$, we have either
\begin{equation}  \label{g'}
\re(\beta)|_{\g t_{\BB R}(X)} \equiv 0
\qquad \text{or} \qquad \re(\beta(X)) \ne 0.
\end{equation}
\end{enumerate}
Clearly, $\g g_{\BB R}'$ is an open subset of $\g g_{\BB R}$;
since $M$ is compact and $\Delta(X)$ is finite, the complement of
$\g g_{\BB R}'$ in $\g g_{\BB R}$ has measure zero.

From now on we will assume that the element $X \in \g g_{\BB R}$ is
not only regular semisimple, but also lies in $\g g_{\BB R}'$.
Then applying the integral localization formula of \cite{L3} to the
integral character formula (\ref{intcharformula1}), with the global formula
for the multiplicities, we can rewrite the formula (\ref{m1}) as
\begin{equation}  \label{m2}
m_{\g b_k}(X) = \chi \bigl( {\cal B}, \bigl( \BB D_{\cal B}
(Rp_*(\widetilde{\cal S})) \bigr)_{O_k} \bigr)
= \chi \bigl( {\cal B},
(j_{O_k \hookrightarrow {\cal B}})_! \circ
(j_{O_k \hookrightarrow {\cal B}})^* \bigl( \BB D_{\cal B}
(Rp_*(\widetilde{\cal S})) \bigr) \bigr).
\end{equation}
Applying Proposition 2.5.11 from \cite{KaScha} to the Cartesian square
$$
\begin{matrix}
O_k \times M & \hookrightarrow & {\cal B} \times M \\
\downarrow & & \downarrow \\
O_k & \hookrightarrow & {\cal B}
\end{matrix}
$$
and using that the projection map $p$ is proper we obtain:
\begin{multline}  \label{m3}
m_{\g b_k}(X) =
\chi \bigl( {\cal B},
(j_{O_k \hookrightarrow {\cal B}})_! \circ
(j_{O_k \hookrightarrow {\cal B}})^* \bigl( Rp_!
(\BB D_{{\cal B} \times M} (\widetilde{\cal S})) \bigr) \bigr)  \\
= \chi \bigl( {\cal B},
(j_{O_k \hookrightarrow {\cal B}})_! \circ Rp_! \circ
(j_{O_k \times M \hookrightarrow {\cal B} \times M})^*
(\BB D_{{\cal B} \times M} (\widetilde{\cal S})) \bigr)  \\
= \chi \bigl( {\cal B}, Rp_! \circ
(j_{O_k \times M \hookrightarrow {\cal B} \times M})_! \circ
(j_{O_k \times M \hookrightarrow {\cal B} \times M})^*
(\BB D_{{\cal B} \times M} (\widetilde{\cal S})) \bigr)  \\
= \chi \bigl( {\cal B} \times M,
(j_{O_k \times M \hookrightarrow {\cal B} \times M})_! \circ
(j_{O_k \times M \hookrightarrow {\cal B} \times M})^*
(\BB D_{{\cal B} \times M} (\widetilde{\cal S})) \bigr)  \\
= \chi \bigl( {\cal B} \times M,
(\BB D_{{\cal B} \times M} (\widetilde{\cal S}))_{O_k \times M} \bigr).
\end{multline}

\separate

Next we will compare this result with the result of application of
the integral localization formula of \cite{L3} to the integral
(\ref{mainintegral}). Let
$$
M = \coprod_{\{x \in M;\: X_M(x) =0 \}} \tilde O_x
$$
be the Bialynicki-Birula decomposition \cite{Bi} of M into attracting sets
(relative to $X$), as described in \cite{L3}.
To obtain this decomposition we pick any
$X' \in \g t_{\mathbb R}(X) \cap \g g_{\mathbb R}'$ in the same connected
component of $\g t_{\mathbb R}(X) \cap \g g_{\mathbb R}'$ and such that
$$
\re \beta(X) >0 \: \Longleftrightarrow \: \re \beta(X') >0
\qquad \text{and} \qquad
\re \beta(X) <0 \: \Longleftrightarrow \: \re \beta(X') <0
$$
for all $\beta \in \Delta(X)$, and the complex 1-dimensional subspace
$\{tX';\: t \in \BB C \} \subset \g g_{\BB C}$ is the Lie algebra of
a closed algebraic subgroup $\BB C^{\times}(X') \subset G_{\BB C}$
isomorphic to $\BB C^{\times}$.
We fix an embedding of $\BB C^{\times}(X') \simeq \BB C^{\times}$ into
$\mathbb C$ so that the tangent map sends
$X' \in T_e (\BB C^{\times}(X'))$ into an element with
nonnegative real part. This embedding allows us to take limits as
$z \in \BB C^{\times}(X')$ approaches to zero,
and for each $x \in M$ with $X_M(x) =0$ we set
$$
\tilde O_x = \{ y \in M;\: \lim_{z \to 0} z^{-1} \cdot y =x \}.
$$
The sets $\tilde O_x$ are smooth locally closed algebraic subvarieties of $M$.

Then
\begin{multline*}
\{x \in {\cal B} \times M;\: X_{{\cal B} \times M} (x) =0 \} \\
= \quad
\{x \in {\cal B};\: X_{\cal B}(x) =0 \} \times \{x \in M;\: X_M(x) =0 \}  \\
= \quad \{\g b_1,\dots,\g b_{|W|}\} \times \{x \in M;\: X_M(x) =0 \}
\end{multline*}
is the set of zeroes of the vector field $X_{{\cal B} \times M}$, and
$$
{\cal B} \times M =
\coprod_{\begin{matrix} k=1,\dots,|W| \\ \{x \in M;\:X_M(x) =0 \} \end{matrix}}
O_k \times \tilde O_x
$$
is the Bialynicki-Birula decomposition of ${\cal B} \times M$
into attracting sets (relative to $X$).
Applying the integral localization formula to the integral (\ref{mainintegral})
yields
$$
\frac {i^n}{(2\pi i)^{d+n}}
\int_{\Lambda} \Bigl( \int_{\g g_{\BB R}} 
\tilde p^* e^{\langle X, \mu_{\cal B}(\zeta) - \lambda_x \rangle +
\sigma_{\cal B} + \tau_{\lambda}}
\wedge \tilde q^* \thom_{T^*M} \wedge \phi(X) \Bigr)_{[2(d+n)]}
= \int_{\g g_{\BB R}} \tilde F_{\theta} \phi
$$
for any $\phi \in \Omega_c^{top}(\g g_{\BB R})$
with $\supp \phi \subset \g g_{\BB R}'$,
where the function $\tilde F_{\theta}$ is
\begin{equation}  \label{Ftilde}
\tilde F_{\theta}(X) =
(-1)^d \sum_{k=1}^{|W|} \sum_{\{x \in M;\:X_M(x) =0 \}} \tilde m_{k,x}(X)
\frac {e^{-\langle X , \lambda_{\g b_k} \rangle}}
{\prod_{\alpha \in \Psi^+_{\g b_k}(X)} \alpha(X)},
\end{equation}
and
\begin{multline*}
\tilde m_{k,x}(X) =
\chi \bigl( {\cal B} \times M, (\BB D_{{\cal B} \times M}
(\widetilde{\cal S}))_{O_k \times \tilde O_x} \bigr)  \\
= \chi({\cal B} \times M, 
(j_{O_k \times \tilde O_x \hookrightarrow {\cal B} \times M})_! \circ
(j_{O_k \times \tilde O_x \hookrightarrow {\cal B} \times M})^*
(\BB D_{{\cal B} \times M} (\widetilde{\cal S})) \bigr).
\end{multline*}
We want to show that $\tilde F_{\theta}(X) = F_{\theta}(X)$. We have:
\begin{multline*}
\sum_{\{x \in M;\:X_M(x) =0 \}} \tilde m_{k,x}(X) =
\sum_{\{x \in M;\:X_M(x) =0 \}}
\chi \bigl( {\cal B} \times M, (\BB D_{{\cal B} \times M}
(\widetilde{\cal S}))_{O_k \times \tilde O_x} \bigr)  \\
= \chi \bigl( {\cal B} \times M,
(\BB D_{{\cal B} \times M} (\widetilde{\cal S}))_{O_k \times M} \bigr)
= m_{\g b_k}(X).
\end{multline*}
Combining this with the equations (\ref{fpcharf}) and (\ref{Ftilde}) we obtain
$$
\tilde F_{\theta}(X) =
(-1)^d \sum_{k=1}^{|W|} m_{x_k}
\frac {e^{-\langle X , \lambda_{\g b_k} \rangle}}
{\prod_{\alpha \in \Psi^+_{\g b_k}(X)} \alpha(X)} = F_{\theta}(X),
$$
which proves (\ref{mainintegral}) at least when
$\supp(\phi) \subset \g g_{\BB R}'$.
Since the function $F_{\theta}$ is known to be a character of some
representation (namely (\ref{vrep})),
it is locally $L^1$ and the main result of \cite{L3} now proves
(\ref{mainintegral}) for all $\phi \in \Omega_c^{top}(\g g_{\BB R})$.

\separate

\begin{section}
{An example on the enhanced flag variety}  \label{example}
\end{section}

The purpose of this section is to show how Theorem \ref{main} can be applied
to smooth quasi-projective varieties which are not necessarily projective.
We let our group $G_{\BB R}$ be
$GL(N, \BB R) \subset GL(N, \BB C) = G_{\BB C}$, for some
$N \in \BB N$; we denote by $B_{\BB C}$ the group of all invertible
lower-triangular matrices and by $N_{\BB C} \subset B_{\BB C}$ the group
of lower-triangular unipotent matrices, we set
$B_{\BB R} = B_{\BB C} \cap G_{\BB R}$ and
$N_{\BB R} = N_{\BB C} \cap G_{\BB R}$.
We will be interested in the complex homogeneous space
$Z_{\BB C} = G_{\BB C} / N_{\BB C}$ and its real submanifold
$Z_{\BB R} = G_{\BB R} / N_{\BB R} \subset Z_{\BB C}$.
The space $Z_{\BB C}$ can be regarded as a fiber bundle over the flag variety
${\cal B} = G_{\BB C} / B_{\BB C}$ with fibers isomorphic to
$H_{\BB C} =_{def} B_{\BB C} /N_{\BB C} \simeq (\BB C^{\times})^N$.
It is easy to see that these fibers can be identified with each other in a
canonical way. This allows us to think of $Z_{\BB C}$ as a trivial principal
$H_{\BB C}$-bundle which is also a $G_{\BB C}$-equivariant fiber
bundle, and the actions of $H_{\BB C}$ and $G_{\BB C}$ commute.
W.~Schmid and K.~Vilonen call $Z_{\mathbb C}$ the {\em enhanced flag variety}
of $\g g_{\BB C}$ \cite{SchV2}.

Let ${\cal S}' = (j_{Z_{\BB R} \hookrightarrow Z_{\BB C}})_*
(\BB C_{Z_{\BB R}})$ denote the sheaf on $Z_{\BB C}$ which is the direct
image of the constant sheaf on $Z_{\BB R}$.
Fix a $\lambda \in \g h_{\mathbb C}^*$ and
denote by $\ker \chi_{\lambda-\rho} \subset {\cal Z}(\g g_{\mathbb C})$
the kernel of the character
$\chi_{\lambda-\rho} : {\cal Z}(\g g_{\mathbb C}) \to \mathbb C$, and take a
${\cal D}_{Z_{\BB C}}$-module
$$
\g M' = {\cal D}_{Z_{\BB C}} /
\gamma_{Z_{\BB C}} (\ker \chi_{\lambda-\rho}) {\cal D}_{Z_{\BB C}}.
$$
Notice that by construction ${\cal Z}(\g g_{\mathbb C})$ acts on $\g M'$
by character $\chi_{\lambda-\rho}$.
We form a virtual representation of $G_{\BB R}$
\begin{equation}  \label{Whittaker}
\rhom^{top}_{{\cal D}_{Z_{\BB C}}} ( \g M' \otimes {\cal S}',
{\cal O}_{Z_{\BB C}^{an}}).
\end{equation}
The composition factors of this representation are known to be
the principal series representations of $G_{\BB R}$ which have the same
character on $\g g_{\BB R}$.
We will obtain an integral formula for the character of this representation.
Of course, the homogeneous space $Z_{\BB C}$ is not compact,
and in order to apply Theorem \ref{main} we need to compactify
$Z_{\BB C}$ first.

Recall that $Z_{\BB C}$ is a principal ($H_{\BB C}, G_{\BB C}$)-equivariant
fiber bundle with commuting $H_{\BB C}$- and $G_{\BB C}$-actions.
Hence we can embed
$H_{\BB C} =_{def} B_{\BB C} /N_{\BB C} \simeq (\BB C^{\times})^N$
into the vector space $\BB C^N$ and
form a $G_{\BB C}$-equivariant vector bundle
$\tilde Z_{\BB C} \simeq \BB C^N \times {\cal B}$
with an embedding $Z_{\BB C} \hookrightarrow \tilde Z_{\BB C}$
of $G_{\BB C}$-equivariant fiber bundles over ${\cal B}$.
Next we form a line bundle $L_0 = \BB C \times {\cal B}$ over ${\cal B}$
with trivial $G_{\BB C}$-action on the first factor and form a
projectivization $M =_{def} \BB P(\tilde Z_{\BB C} \oplus L_0)$.
Then $M$ is a smooth complex projective variety, $G_{\BB C}$ acts on $M$
algebraically, and $M$ contains $Z_{\BB C}$ as a dense open subset.
Define
$$
\g M = {\cal D}_M / \gamma_M (\ker \chi_{\lambda-\rho}) {\cal D}_M
\qquad \text{and} \qquad
{\cal S} = (j_{Z_{\BB C} \hookrightarrow M})_! {\cal S}'
= (j_{Z_{\BB R} \hookrightarrow M})_! (\BB C_{Z_{\BB R}}).
$$

For convenience we restate Theorem 5.12 of \cite{KaSchm}:
\begin{thm}
Let $f: X \to Y$ be a $G_{\BB C}$-equivariant morphism between complex
algebraic, quasi-projective $G_{\BB C}$-manifolds $X$, $Y$.
If $f$ is projective, there exists an isomorphism
$$
\rhom^{top}_{{\cal D}_X} ( \g M \otimes f_{an}^{-1} {\cal T},
{\cal O}_{X^{an}}) \quad \simeq \quad
\rhom^{top}_{{\cal D}_Y} \Bigl ( \int_f \g M \otimes {\cal T},
{\cal O}_{Y^{an}} \Bigr ) [-\dim X/Y],
$$
functorially in ${\cal T} \in D^b_{G_{\BB R}, \BB R-c}(\BB C_Y)$ and
$\g M \in D^b_{G_{\BB C}, coh}({\cal D}_X)$
(the bounded derived category of $\modg^{coh}({\cal D}_X)$).
If $f$ is smooth,

$$
\rhom^{top}_{{\cal D}_X} (Lf^* \g M \otimes {\cal T}, {\cal O}_{X^{an}})
\quad \simeq \quad
\rhom^{top}_{{\cal D}_Y} (\g M \otimes R(f_{an})_! {\cal T},
{\cal O}_{Y^{an}}) [-2\dim X/Y],
$$
functorially in $\g M \in D^b_{G_{\BB C}, coh}({\cal D}_Y)$
and ${\cal T} \in D^b_{G_{\BB R}, \BB R-c}(\BB C_X)$.
\end{thm}

Applying the second part of this theorem to the open inclusion
$Z_{\BB C} \hookrightarrow M$
we can rewrite our virtual representation (\ref{Whittaker}) as
$$
\rhom^{top}_{{\cal D}_{Z_{\BB C}}} ( \g M' \otimes {\cal S}',
{\cal O}_{Z_{\BB C}^{an}}) \quad \simeq \quad
\rhom^{top}_{{\cal D}_M} ( \g M \otimes {\cal S}, {\cal O}_{M^{an}}).
$$
The ${\cal D}_M$-module $\g M$ is ${\cal Z}(\g g_{\mathbb C})$-finite and
lies in $\modg^{coh,\lambda-\rho}({\cal D}_M)$ essentially by construction.
It follows from the Bruhat decomposition of $G_{\mathbb C}$ that
$B_{\BB C}$ acts on $M$ with finitely many orbits; this implies that
$\g M$ is admissible.
It is easy to see that each maximal torus $T_{\BB C}$ acts on $M$ with
finitely many fixed points.
Therefore Theorem \ref{main} applies here and we obtain an integral
formula for the character of the virtual representation (\ref{Whittaker}).
This is probably the most complicated formula for the principal series
character there is.

\separate

\end{document}